\documentclass[12pt,a4paper]{article}
\usepackage[english]{babel}
\usepackage{verbatim}
\usepackage{bm, color}

\usepackage{amsmath}
\usepackage{amssymb}
\usepackage{amsfonts}
\usepackage{amsbsy}

\usepackage{graphicx}
\usepackage{graphics}
\usepackage{subfigure}
\usepackage{epsfig}

\usepackage{array}
\usepackage{booktabs}
\usepackage{caption}
\usepackage{tabularx}

\newcommand{\R}{{\mathbb R}}

\begin{document}

\title{A fully  semi-Lagrangian discretization\\
for the 2D Navier--Stokes equations\\
in the vorticity--streamfunction formulation}

\author{Luca Bonaventura$^{(1)}$,\ \  Roberto Ferretti$^{(2)}$,\\ Lorenzo Rocchi$^{(2)}$}

\maketitle

\begin{center}
{\small
$^{(1)}$ MOX -- Modelling and Scientific Computing, \\
Dipartimento di Matematica, Politecnico di Milano \\
Via Bonardi 9, 20133 Milano, Italy\\
{\tt luca.bonaventura@polimi.it}
}
\end{center}

\begin{center}
{\small
$^{(2)}$
Dipartimento di Matematica e Fisica\\
Universit\`a degli Studi Roma Tre\\
 L.go S. Leonardo Murialdo 1,
00146, Roma, Italy\\
{\tt ferretti@mat.uniroma3.it, lor\_89@fastwebnet.it}
}
\end{center}

\date{}


\vspace*{0.5cm}


\vspace*{0.5cm}

\pagebreak

\abstract{A numerical  method for the two-dimensional, incompressible Navier--Stokes
equations in vorticity--streamfunction form is proposed, which employs semi-Lagrangian discretizations for both the advection and diffusion terms, thus achieving unconditional stability without
the need to solve linear systems beyond that required by the Poisson solver for the reconstruction of the streamfunction. A description of the  discretization of Dirichlet boundary conditions
for the semi-Lagrangian approach to diffusion terms is  also presented.
 Numerical experiments on classical benchmarks for incompressible flow in simple geometries
  validate the proposed method.} 

\pagebreak

\section{Introduction}
\label{intro} \indent
Over the last 30 years, Semi-Lagrangian (SL) schemes  have been extremely successful for advection dominated problems,  allowing for example major reductions in the computational cost of operational
weather predictions, see e.g. \cite{bonaventura:2012}, \cite{staniforth:1991}. The potential advantage of SL schemes comes from their unconditional stability with respect to the Courant number, which is in turn achieved by exploiting  approximations of  the characteristics of the advection equation.

Several extensions of SL schemes to advection--diffusion problems have been proposed in the last two decades.
Some of these techniques are based on a splitting of the evolution operator, in which only the advective part is treated in SL form, whereas others treat the diffusion term in a consistently SL form by using multiple characteristics, as it will be shown in Section~\ref{fullysl}. A review of these proposals is presented, e.g., in~\cite{falcone:2013}. 

Within the second line of work, which has been traditionally motivated by stochastic arguments, we quote here the seminal paper \cite{camilli:1995}, along with the more recent contributions in \cite{milstein:2002}, \cite{milstein:2000}, \cite{milstein:2001} and \cite{ferretti:2010}. Due to their stochastic origin, all these works treat diffusion operators in trace form, but  extensions to divergence form operators has also been proposed in \cite{bonaventura:2014}, \cite{bonaventura:2016a}. In these papers, SL methods have proved to be accurate and efficient for linear and nonlinear
advection--diffusion problems.

In the present work we are interested in applying the same approach to the incompressible Navier--Stokes (NS) equations
\begin{eqnarray}\label{NSE}
 &&\frac{\partial\bm{u}}{\partial t}+\left(\bm{u}\cdot\nabla\right)\bm{u}-\nu\Delta\bm{u}+\nabla p=0, \nonumber \\
&&\nabla\cdot \bm u = 0,
\end{eqnarray}
on a bounded domain $\Omega\subset\R^2$, with proper initial and boundary conditions. 
SL techniques of both the forms outlined above have been proposed for this problem, we refer 
here for example to the papers \cite{maday:1990}, \cite{xiu:2001}, \cite{xu:2002}, \cite{xiu:2005} 
and the review in \cite{bermejo:2016}, in which the application of the SL method is restricted to the advective part, and  to \cite{belopolskaya:2003}, \cite{milstein:2012}, which present and analyze a fully SL approach based on a stochastic framework. In particular, these latter works present a theoretical analysis for the time-discrete scheme, although with  limited numerical validation.

In the present work, we study the application of a fully SL scheme, much in the same spirit of \cite{belopolskaya:2003}, \cite{milstein:2012}, to the two-dimensional NS equations from a more numerical perspective. We  will use the vorticity and streamfunction formulation for simplicity,
as a first step towards the application of the same approach  in the context of a projection method for the NS equations in their standard velocity-pressure formulation.
We will show that the fully SL approach yields an explicit discretization of both advection and diffusion terms
 with very mild stability restrictions,  that  can achieve higher order spatial accuracy in a practically and conceptually simple way, while reducing the computational cost of the advection--diffusion step. 
 The scheme will be constructed and its consistency analyzed. Moreover, we will provide details on the implementation of boundary conditions, which was not discussed in detail in our 
 previous papers, and perform some classical numerical tests to validate the method.
Numerical results  show that the method yields good quantitative agreement with reference
numerical solution of classical benchmarks, while allowing the use of time steps several times
larger than those of a standard explicit scheme.

The outline of the paper is the following. Section \ref{vortstream} recalls the vorticity--streamfunction formulation of 2D NS equations, along with the related boundary conditions. Section \ref{fullysl} describes the SL 
advection--diffusion solver, along with its consistency analysis and the  implementation of the boundary conditions.
A numerical validation of the proposed approach on classical benchmarks 
is presented in Section \ref{tests}, while some conclusions and perspectives for future developments
are outlined in Section \ref{conclu}.

\section{The Navier--Stokes equations in the vorticity--streamfunction formulation}
\label{vortstream} \indent
 
We briefly recall here the vorticity--streamfunction (VS) formulation of the NS equations.
The basic idea, see \cite{chorin:1993} for an in-depth discussion, is that any two-dimensional divergence-free vector field $\bm u=(u_1, u_2)$ can be obtained as
\begin{equation}\label{grad_ort}
u_1=\frac{\partial\psi}{\partial x_2},\qquad u_2=-\frac{\partial\psi}{\partial x_1},
\end{equation}
for a suitable streamfunction $\psi$. Now, since $\bm u$ is two-dimensional, its vorticity is a scalar given by
\begin{equation}
\omega=\frac{\partial u_2}{\partial x_1}-\frac{\partial u_1}{\partial x_2},
\end{equation}
which, together with \eqref{grad_ort}, gives
\begin{equation}\label{poisson}
-\Delta\psi=\omega.
\end{equation}
On the other hand, by taking the curl of both sides in \eqref{NSE} we get
\begin{eqnarray}\label{dif_tra}
0 & = & \frac{\partial\omega}{\partial t}-\nu\Delta\omega
+\frac{\partial\psi}{\partial x_2}\frac{\partial\omega}{\partial x_1}
-\frac{\partial\psi}{\partial x_1}\frac{\partial\omega}{\partial x_2} \nonumber \\
& = & \frac{\partial\omega}{\partial t}-\nu\Delta\omega
+u_1\frac{\partial\omega}{\partial x_1}+u_2\frac{\partial\omega}{\partial x_2},
\end{eqnarray}
i.e., an advection--diffusion equation for the vorticity. Then, the VS formulation combines \eqref{grad_ort}, \eqref{poisson} and \eqref{dif_tra} in the form
\begin{eqnarray}\label{VS}
&&\frac{\partial\omega}{\partial t}+\bm{u}\cdot\nabla\omega-\nu\Delta\omega=0\nonumber\\
&& -\Delta\psi=\omega \ \ \ \ \ \  \bm{u}=\nabla^{\perp}\psi
\end{eqnarray}
in which the symbol $\nabla^\perp$ denotes in compact form the operator defined by \eqref{grad_ort}. Note that the incompressibility condition is ensured by \eqref{grad_ort} and that the advantage of treating a scalar problem for the vorticity only occurs  in two space dimensions.  

The precise derivation of boundary conditions for the vorticity--streamfunction formulation is not trivial,
see e.g. the discussion in \cite{e:1996}, \cite{quartapelle:2013}. Here, we will simply note that
no slip  boundary conditions amount to assigning

\begin{equation}\label{VS_bc}
\psi=0, \ \ \ \ \ \ \  \nabla^{\perp}\psi \cdot {\bf n}=0,
\end{equation}
where ${\bf n} $ denotes the outer normal at the boundary,   while more general non homogeneous 
boundary conditions are assigned by setting non zero values for the normal derivative of the streamfunction.

 \section{The fully semi-Lagrangian numerical method}
 \label{fullysl} \indent
 
 Our aim in this work is to
 employ the fully SL advection--diffusion method described in \cite{bonaventura:2014}
  in the framework a  discretization of the NS equations in VS formulation. The general idea underlying this technique is to  approximate the diffusion term by a convex combination of the known values at the locations of a time dependent stencil. Similar ideas, often justified by probabilistic arguments, have been proposed in various independent contexts,  see \cite{camilli:1995},   \cite{ferretti:2010}, \cite{milstein:2013},  \cite{teixeira:1999}.

Given the advection--diffusion equation
\begin{eqnarray}\label{trasp}
&&      \omega_t+\bm u\cdot \nabla\omega-\nu\Delta\omega = 0,\ \ \  (x,t)\in\R^2 \times [0,T] \nonumber \\
 &&     \omega(x,0)=\omega_0(x) \ \ \  x\in\R^2,
\end{eqnarray}
we first consider the inviscid case $\nu=0$. The construction of large time-step schemes for \eqref{trasp} stems from the application of the classical  method of characteristics. The system of characteristic curves $X(x,t;s)$ for \eqref{trasp} is defined by the solutions of:
\begin{eqnarray}\label{caratt}
&&\displaystyle\frac{d}{ds}X(x,t;s)=\bm u(X(x,t;s),s),\nonumber\\
&&X(x,t;t)=x.
\end{eqnarray}
The solution of \eqref{trasp} is constant along such curves, which means that it satisfies the relationship
\begin{equation}\label{upw}
\omega(x,t)=\omega(X(x,t;t-\Delta t),t-\Delta t).
\end{equation}
Discretizing the representation formula \eqref{upw} we obtain the advective SL approximation. More precisely, we denote by $\Delta t$ and $\Delta x$ respectively the time and space discretization steps, with $t_n=n\Delta t$ for $n\in [0,T/\Delta t]$, and a space grid of points $x_i$. The characteristics $X$ defined by \eqref{caratt} will be replaced by their numerical approximations $X^\Delta. $ We will also use the shorthand notation
\begin{equation}\label{z_i}
z_i^{n+1}=X^\Delta(x_i,t^{n+1};t^n)
\end{equation}
to denote the foot of the approximate characteristic starting from $(x_i,t_{n+1}).$
In this work, we have employed characteristics approximations based on the standard
explicit Euler and Heun schemes. These methods are applied with substepping, as described e.g. in \cite{casulli:1990}, \cite{giraldo:1999}, \cite{rosatti:2005}. More precisely, for the 
approximation of characteristics  a time step $\Delta \tau $ is employed,
that is a fraction of the total time step $\Delta t.$  This time step is required to comply with a local CFL restriction.
As discussed in \cite{rosatti:2005}, the value of 
$\Delta \tau $ can vary along each approximated characteristic, so that $\Delta \tau << \Delta t $  only in those
(usually small) regions with really large Courant numbers. During these substeps, the velocity field
$\bm u $ is frozen at a constant time level for simplicity. Since employing $\bm u^n $  
entails at most first order accuracy, in the case of the Heun scheme the extrapolation
$$
\bm u^{n+\frac 12}\approx\frac 32 \bm u^n - \frac 12 \bm u^n
$$
is employed, as common in the meteorological literature, see e.g. \cite{staniforth:1991}.

 For inviscid problems with $\nu= 0$, \eqref{upw} is discretized by replacing the exact upwinding $X$ with $X^\Delta$, and the value of $u$ at the foot of a characteristic with an interpolation $I_p$:
\begin{equation}\label{sl}
\omega_i^{n+1}=I_p[\Omega^n] \left(z_i^{n+1}\right)
\end{equation}
where $\omega_i^{n+1}$ is the approximation of $\omega(x_i,t^{n+1})$, $\Omega^n$ is the vector of node values $\omega_i^n$, and $I_p$ is an interpolation operator (e.g., a polynomial interpolation of degree $p$) which is assumed to satisfy the condition
\[
I_p[\Omega](x_i) = \omega_i.
\]
As customary in the implementation of semi-Lagrangian schemes, whenever
$z_i^{n+1}$ falls outside of the computational domain, it is redefined as the closest boundary
point. Even though this implies in principle an error that is first order in time, its impact is minimized
by the use of substepping approaches for the approximation of the trajectories, see e.g. the discussion
in \cite{rosatti:2005}.

When $\nu\neq 0$, the scheme can be modified along the lines proposed
in \cite{bonaventura:2014}, \cite{ferretti:2010}, in order to introduce an approximation of the diffusion
term. In fact, a first-order consistent discretization of the term $\nu\Delta\omega$ is obtained by replacing the interpolation of the numerical solution at $z_i^{n+1}$ with an average of interpolated values obtained adding to $z_i^{n+1}$ a second displacement of the form
\[
\delta_k = \pm \delta e_j, \quad (k=1,\ldots,4)
\]
for all combination of both the sign and the index $j=1,2$, and with
\[
\delta=\sqrt{4\nu\Delta t}.
\]  
Notice that, in an alternative form of the diffusive part of the scheme, 
the displacements $\delta_k$ could be rather defined as
\[
\delta_k =\frac{\delta}{\sqrt 2} \begin{pmatrix} \pm 1 \\ \pm 1\end{pmatrix},
\]
which corresponds to the definition given in \cite{belopolskaya:2003, milstein:2012}. While this definition is perfectly equivalent in two space dimensions (the stencil of points is simply rotated), in three space dimensions this would require a higher number of interpolations at each node (eight interpolation points instead of six).

\begin{figure}
\centering
\includegraphics[height=4cm]{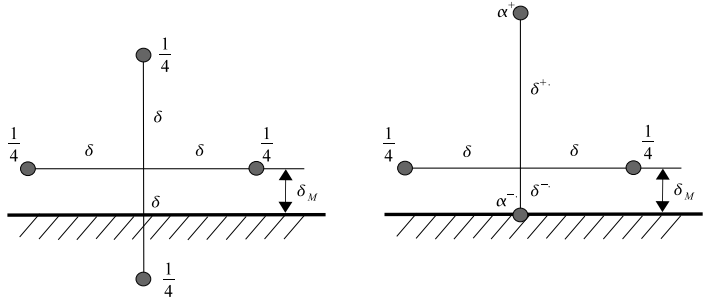}
\caption{Modification of weights and displacements near the boundary}\label{fig:BoundStencil}
\end{figure}

The resulting scheme then reads 
\begin{equation}\label{sl_dt}
\omega_i^{n+1} = \frac{1}{4} \sum_{k=1}^4 I_p[\Omega^n] \left(z_i^{n+1}+\delta_k \right),
\end{equation}
in which $z_i^{n+1}$ is defined by \eqref{z_i}. 
In case some of the displacements $\delta_k$ identifies a location out of the computational domain,
the displacement is redefined so that the point $z_i^{n+1}+\delta_k$ is the boundary point along the line
connecting $z_i^{n+1}$ and $z_i^{n+1}\pm\delta e_j.$ For these points, the scheme is redefined as
\begin{equation}\label{sl_dt_bound}
\omega_i^{n+1} = \sum_{k=1}^4 \alpha_k I_p[\Omega^n] \left(z_i^{n+1}+\delta_k \right),
\end{equation}
where $\sum_{k=1}^4\alpha_k=1$. In order to define the displacements $\delta_k$ and weights $\alpha_k$ in a consistent form, we consider without loss of generality the simple case depicted in Figure \ref{fig:BoundStencil}. 
A pure diffusion operator  is discretized by the abstract scheme
\begin{eqnarray*}
u(x,t_{n+1}) & \approx & \frac{1}{4} u(x+\delta e_1,t_n) + \frac{1}{4} u(x-\delta e_1,t_n) \\
&& + \alpha^+ u(x+\delta^+e_2,t_n) + \alpha^- u(x-\delta^-e_2,t_n),
\end{eqnarray*}
in which we enforce the constraint
\begin{equation}\label{vincolo}
\delta^-\le \delta_M,
\end{equation}
according to the scheme of Fig. \ref{fig:BoundStencil}.
Following \cite{bonaventura:2014}, consistency is achieved by imposing conditions on the moments
of the discrete mass distribution defined by the weights $\alpha^\pm $ and displacements
$\delta^\pm. $ Requiring the first and second moment
 to coincide with the original ones, we obtain the conditions:
\begin{eqnarray}\label{sistema}
&& \alpha^+ + \alpha^- = 1/2 \nonumber\\
&& \alpha^+\delta^+ - \alpha^-\delta^- = 0 \\
&& \alpha^+({\delta^+})^2 + \alpha^-({\delta^-})^2 = 2 \nu\Delta t \nonumber.
\end{eqnarray}
The first condition is satisfied if
\begin{equation*}
\alpha^- = \frac{1}{2}-\alpha^+.
\end{equation*}
Using the relationship in the second equation of \eqref{sistema}, we obtain
\begin{eqnarray*}
0 & = & \alpha^+\delta^+ - \alpha^-\delta^- \\
& = & \alpha^+\delta^+ - (1/2 - \alpha^+)\delta^-
\end{eqnarray*}
which gives in turn
\begin{equation}\label{pesi}
\alpha^+ = \frac{1}{2}\frac{\delta^-}{\delta^+ + \delta^-}, \quad \alpha^- = \frac{1}{2}\frac{\delta^+}{\delta^+ + \delta^-}.
\end{equation}
Using now \eqref{pesi} in the third equation of \eqref{sistema}, we obtain
$$
2\Delta t \nu  =  \frac{1}{2}\frac{\delta^+ {\delta^-}^2 + \delta^- {\delta^+}^2}{\delta^+ + \delta^-} 
 =  \frac{1}{2}\delta^+ \delta^-,
$$
that is,
\begin{equation}\label{spost}
\delta^+=\frac{4\Delta t \nu}{\delta^-}.
\end{equation}
To sum up, consistency under the constraint \eqref{vincolo} would be obtained by setting
\begin{equation}\label{soluz}
\delta^-=\delta_M, \quad \delta^+=\frac{4\Delta t \nu}{\delta_M}, \quad \alpha^-=\frac{1}{2}\frac{\delta^+}{\delta^++\delta^-}, \quad \alpha^+=\frac{1}{2}-\alpha^-.
\end{equation}
Note that, in general, \eqref{sistema} ensures only consistency of order 1/2 (see the analysis in \cite{bonaventura:2014}). Consistency of first order requires the further condition
\[
\alpha^+{\delta^+}^3 - \alpha^-{\delta^-}^3 = O(\Delta t^2),
\]
which is false in general for the solution \eqref{soluz}. Therefore, we should expect the consistency rate of the diffusive term to drop to $\Delta t^{1/2}$ at points where the constraint \eqref{vincolo} holds.

Denoting by $A$ a numerical approximation of the operator $-\Delta$ (e.g., via a five-point discrete Laplacian), by $D^\perp$ an approximation of the operator $\nabla^\perp$ (e.g., via centered differences), and by $\Psi^{n+1}$ and $\bm U^{n+1}$ the vectors of  the node values $\psi_i^{n+1}$ and $\bm u_i^{n+1},$ respectively,  the final form of the scheme is therefore
\begin{eqnarray}\label{VS_discrete}
&& \omega_i^{n+1} =  \sum_{k=1}^4\alpha_k I_p[\Omega^n] \left(z_i^{n+1}+\delta_k \right), \nonumber \\
&& A\Psi^{n+1}=\Omega^{n+1},\\
&& \bm U^{n+1}=D^\perp\Psi^{n+1}, \nonumber
\end{eqnarray}
with initial conditions given by $\Omega^0 = A\Psi^0, $
$\bm U^0 = D^\perp \Psi^0, $ where $\Psi^0 $ denotes the vector of the initial streamfunction values.

\begin{figure}
\centering
\includegraphics[height=4cm]{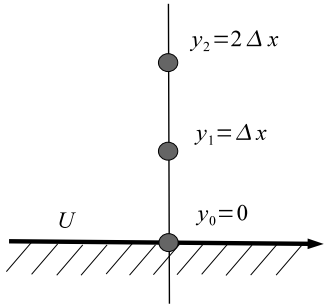}
\caption{Node arrangement for the Thom boundary conditions}\label{fig:thom}
\end{figure}

 For a complete definition of the numerical method, boundary conditions for vorticity
 must be assigned, which do not follow immediately from the boundary conditions
 for the streamfunction. A wide range of possibilities is discussed in \cite{e:1996}. Here,
 we will employ one of the simplest formulations, corresponding to the so called
  Thom boundary conditions \cite{thom:1933}, which are obtained by converting (via Poisson's equation) 
  boundary conditions for the derivative of $\psi$ into Dirichlet boundary conditions for $\omega$.
More specifically,
for a linear portion of the boundary on which a tangential (possibly zero) speed $U$ is imposed,  the Thom boundary conditions read
\begin{equation}\label{thom}
\omega_0 = -\frac{2}{\Delta y^2}(\psi_1 - \psi_0) + \frac{2U}{\Delta y}.
\end{equation}
A sketch of the reference geometry is given in Fig.~\ref{fig:thom}. Note that the Thom conditions are directly written in terms of the approximate solution on the grid, in which the subscript denotes the index of the node with respect to the boundary, so that $\omega_0$, $\psi_0$ refer to the boundary and $\psi_1$ to the first internal layer of nodes.
Note also that, in \eqref{VS_discrete}, the computation of \eqref{thom} should be performed at the start of time step $n+1$ after the Poisson equation has been solved for $\Psi^n$ at time step $n$ (or from the initial data if $n=0$). The SL solution of the advection--diffusion equation should then be computed only at internal nodes, using the values provided by \eqref{thom} as boundary data.  

 \section{Consistency analysis of the method}
 
We give now a simplified consistency analysis for the global method, assuming that it is posed on the whole of $\R^2$ with a uniform orthogonal grid. The analysis of the advection--diffusion step has been carried out elsewhere (see \cite{ferretti:2010, bonaventura:2014}), and shows that the error introduced between two successive time steps is
\[
\epsilon_1 \le C\left(\Delta t^2+ \Delta x^{p+1}\right).
\]
Using a five-point laplacian, the Poisson equation is solved with $O(\Delta x^2)$ accuracy. On the other hand, since the right-hand side is perturbed by an error $\epsilon_1$, the error introduced on the streamfunction is
\[
\epsilon_2 \le C\left(\Delta t^2 + \Delta x^{p+1} + \Delta x^2 \right).
\]
Finally, the error introduced by the operator $D^\perp$, when implemented by centered differences, is again of order $O(\Delta x^2)$. Since for the order of interpolation we have $p\ge 1$, we can drop the term $\Delta x^{p+1}$, and dividing by $\Delta t$ we finally obtain the local truncation error
\[
\tau(\Delta x, \Delta t) \le C\left(\Delta t+ \frac{\Delta x^2}{\Delta t}\right),
\]
which results in a first-order scheme when working at constant Courant number. A comparison with more conventional advection--diffusion solvers (see again \cite{bonaventura:2014}) shows that, despite being formally low-order, this strategy provides a good absolute accuracy. In addition, as usual in SL schemes, 
no stability restrictions on the Courant number are introduced.
However, since the numerical stencil consists of points $O(\sqrt{\Delta t})$ apart, some care should be taken to avoid under-resolution of the smaller spatial scales. The analysis carried out in \cite{ferretti:2010} shows that, in order to have a smooth numerical domain of dependence at time $T$, the compatibility condition
\[
\frac{\nu^{1/3}\Delta t}{T^{2/3}\Delta x^{2/3}} \ll 1
\]
should be enforced on the discretization steps. In particular, this implies
that time steps with Courant numbers quite larger than one are still acceptable,
so that the efficiency of the semi-Lagrangian scheme for the advection term is not lost.

 \section{Numerical results}
 \label{tests} \indent
 
 The method outlined in the previous sections has been implemented so far in a rather elementay way,
 employing structured Cartesian meshes on a rectangular domain.
 The boundary conditions described in section \ref{fullysl} were employed.
  In one numerical experiments,
 the mesh employed had constant step in either direction. The Poisson equation 
 was discretized in this case  via a standard five-point finite difference scheme, and the
 operator $D^\perp$   by centered finite differences. For the interpolation operator
 $I_p,$ the monotonized bicubic interpolant implemented in the MATLAB native function {\tt interp2}
 was employed. Also piecewise linear interpolation was tested, but this resulted in an overly dissipative
 method, as well known in the SL literature for the advective terms (see e.g. \cite{staniforth:1991}) and already
 shown for the diffusive terms in \cite{bonaventura:2014}. The corresponding results are not reported
 in the following.
In the other tests, meshes with non constant steps were employed, in order to achieve higher resolution at the
boundaries, where vorticity production takes place. In this case, first order finite differences 
approximations were employed both for the Poisson equation and the $D^\perp$ operator. For both the constant
and non constant step cases, the finite difference approximations employed for the Poisson equation
yield matrices that coincide with those that would arise from a linear finite element discretization
on the same mesh. Furthermore, for the   meshes with non constant steps, the 
the cubic spline option in {\tt interp2} was employed.

\subsection{Test case with analytic solution}
We consider the test case discussed e.g. in \cite{liu:2000}. In this test, that was run
with $\nu=0.02,  $ over the time interval $t\in[0,4],  $ the computational
domain is given by $[0,2\pi]\times[0,2\pi], $ with periodic boundary conditions.
The exact solution is given by
$  \omega(x,y,t)=\sin{x}\sin{y}\exp{(-2\nu t)}. $
For this test, Cartesian meshes with  constant node spacing were employed, as described above.
The first order convergence behaviour is confirmed, but it can be seen that time steps up to four times
those of explicit schemes can be used without any loss in accuracy.

\begin{table}
\label{anaytic_res1}
    \begin{tabular}{|c|c|c|c|c|c|}\hline
Nodes &  $ u\Delta t/\Delta x$  &   $ \nu\Delta t/(2\Delta x^2)$  & $L^\infty $ rel. err. &   $L^2 $ rel. err. &$ p_{2}$ \\
\hline 
$50\times50$  & $2.6$ & $1$ &$1.12\cdot 10^{-2}$&  $1.68\cdot 10^{-2}$ & - \\ 
$100\times100$ & $2.6$&$2$ &$5.44 \cdot 10^{-3}$ & $7.54 \cdot 10^{-3}$&$1.15 $ \\ 
$200\times200$ & $2.6$&$4$ &$2.58 \cdot 10^{-3}$& $3.57 \cdot 10^{-3}$ &$1.08$ \\ 
\hline
\end{tabular}
\caption{Convergence test in benchmark with analytic solution}
\end{table}

\subsection{Lid driven cavity}
This classical benchmark for numerical approximations of the Navier Stokes equations
has been discussed in detail for the vorticity and streamfunction formulation in \cite{e:1996}.
Reference solutions obtained with high order spectral methods and other accurate
techniques are reported, among many others, in
\cite{botella:1998}, \cite{bruneau:2006}. In our assessment, we focus on low or moderate values
of the Reynolds number. Indeed, it is known that for $Re\approx 8000$ the flow becomes
turbulent and unpredictable for this test, so that a more sophisticated statistical assessment would have
to be carried out. Furthermore, the proposed method is expected to yield significant efficiency gains
especially in  laminar regimes, in which the stability restrictions due to the viscous term can
affect the efficiency of standard explicit methods.

\begin{figure}
\centering
  \includegraphics[width=0.45\textwidth]{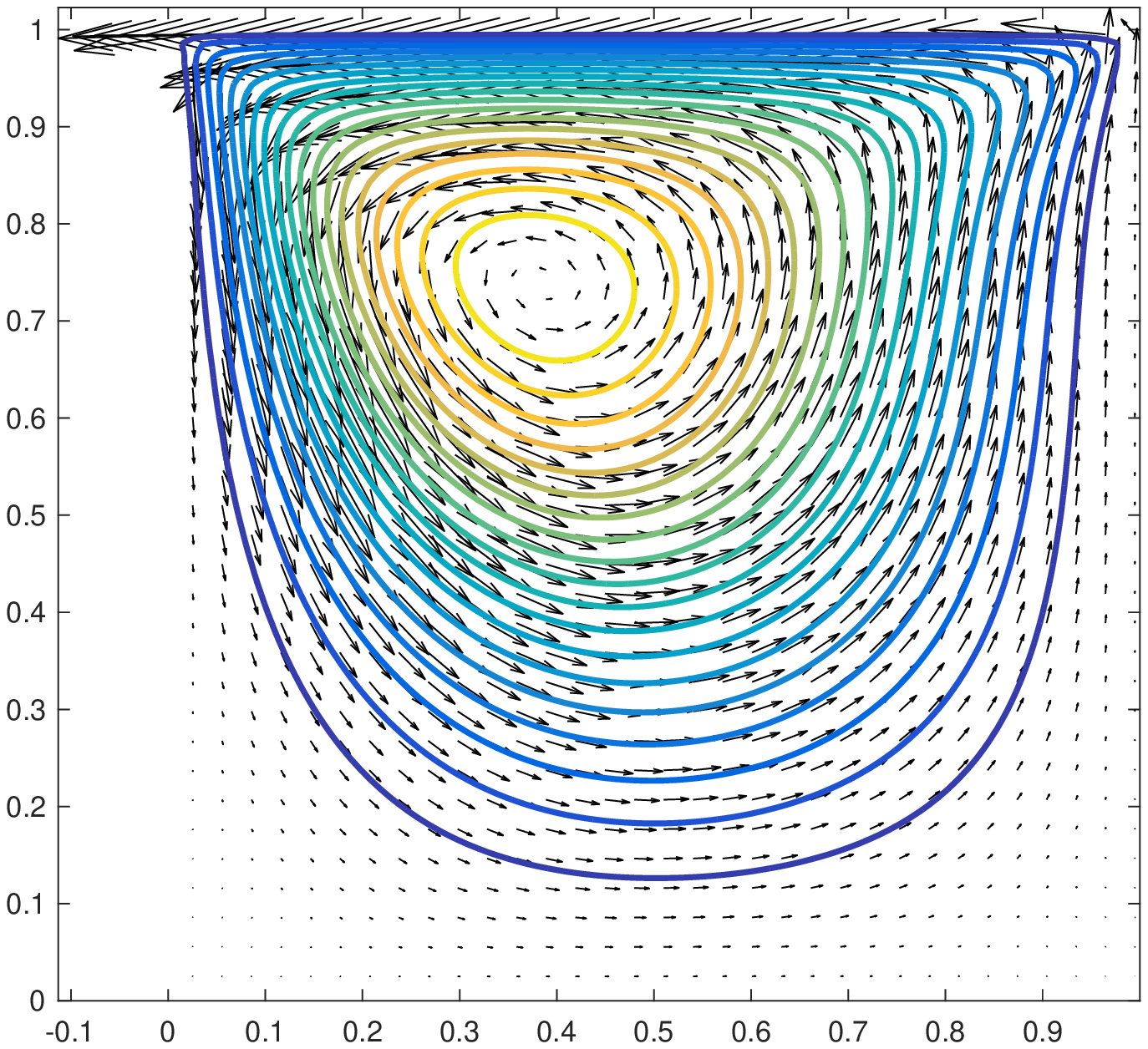}
  \includegraphics[width=0.45\textwidth]{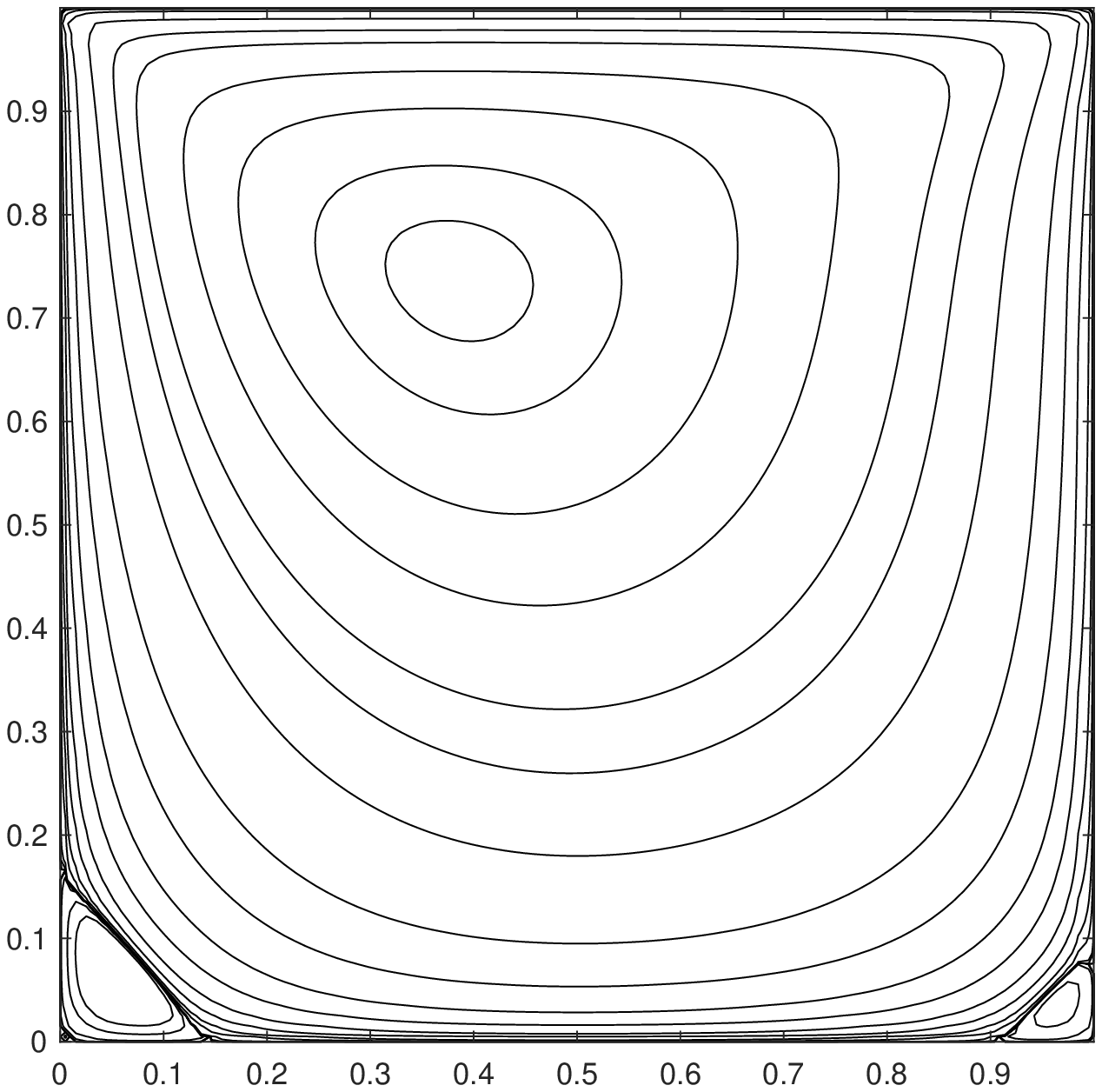}
\caption{Lid driven cavity benchmark at   $Re=100:$ a) flow field b) streamfunction.}
\label{fig:flow_re100}
\end{figure}

An example of the flow field and streamfunction obtained in the case $Re=100$ at steady state
 (approximately
$T=100$) is shown in Figure \ref{fig:flow_re100}. The streamfunction is plotted with the same isoline
values  as suggested in \cite{bruneau:2006}. The results were obtained
using 100 nodes in each direction and a mesh 
whose 2 nodes closest to the boundary are separated by a spacing smaller than that used in the
rest of the domain.
The time step employed was such that $ \nu\Delta t/(2\Delta x^2)\approx 4,  $ 
 $ \|\mathbf{u}\|\Delta t/\Delta x\approx 1 $
for the smallest  $\Delta x $ values and the solution was computed until steady state was achieved
up to a tolerance of $10^{-7}.$
At a more quantitative level, the maximum horizontal velocity value along the centerline of the cavity
was $u_{max}=0.21458,$ while the maximum and minimum vertical velocity values along the centerline of the cavity
were $v_{max}=0.17534,$ $v_{min}=-0.24613,$  respectively. These values imply a relative
error of order $10^{-3}$ with respect to the reference steady state solution in \cite{botella:1998}.
The vorticity value at the center of the cavity was computed as  $\omega_{cen}=1.13370,$ 
which implies a relative
error of order $10^{-2}$ with respect to the same reference  solution.

\begin{figure}
\centering
  \includegraphics[width=0.45\textwidth]{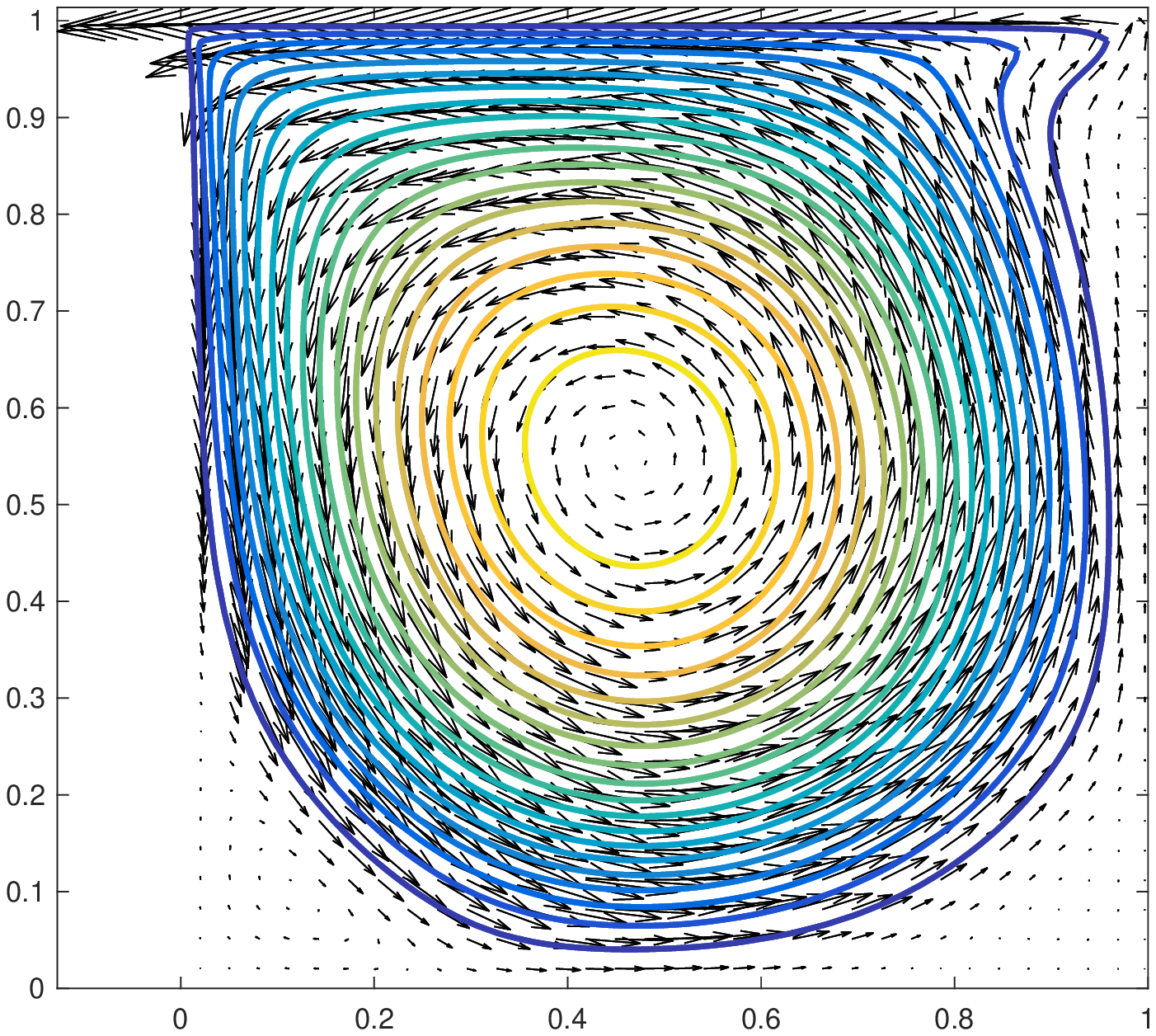}
  \includegraphics[width=0.45\textwidth]{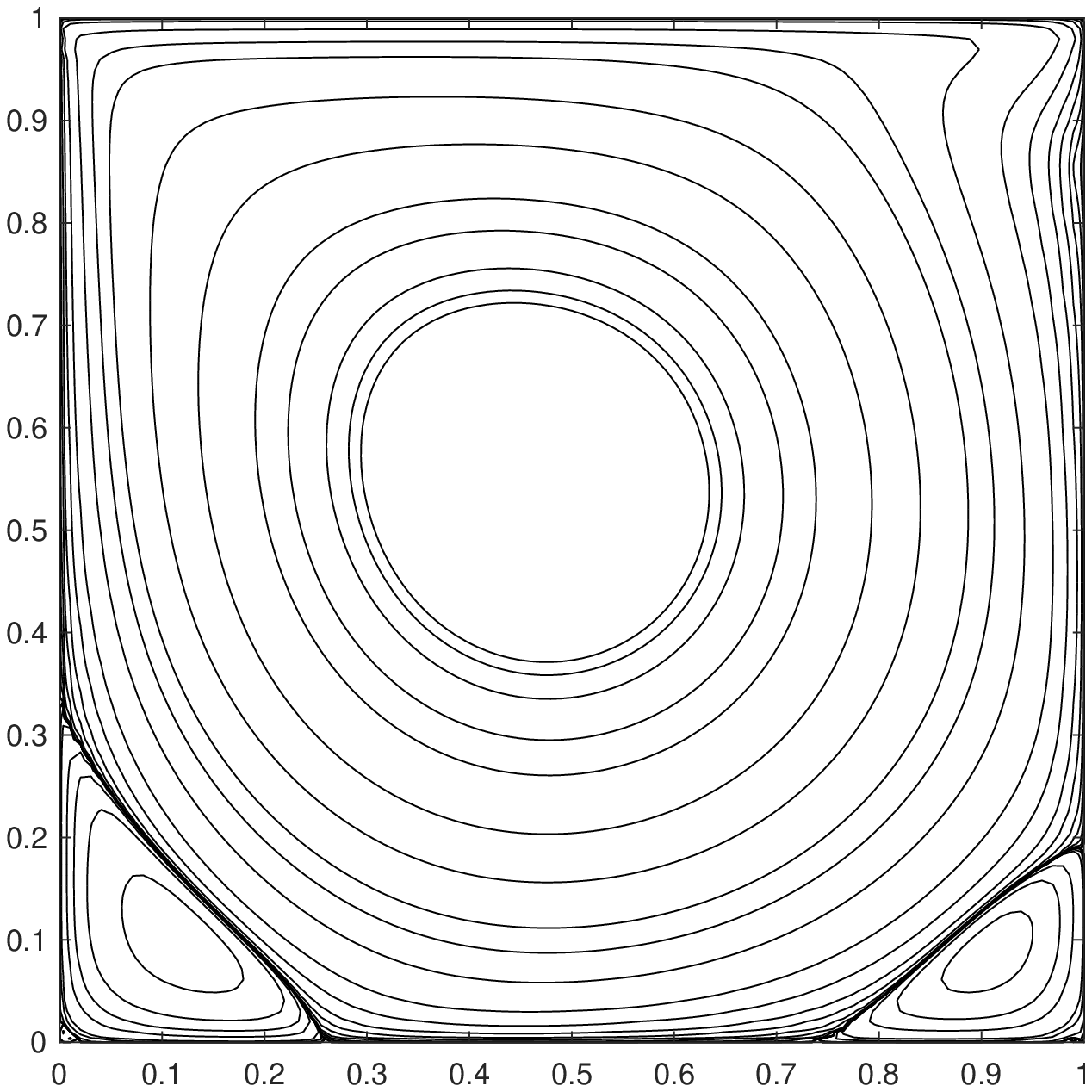}
\caption{Lid driven cavity benchmark at   $Re=1000:$ a) flow field b) streamfunction.}
\label{fig:flow_re1000}
\end{figure}

\begin{figure}
\centering
  \includegraphics[width=0.45\textwidth]{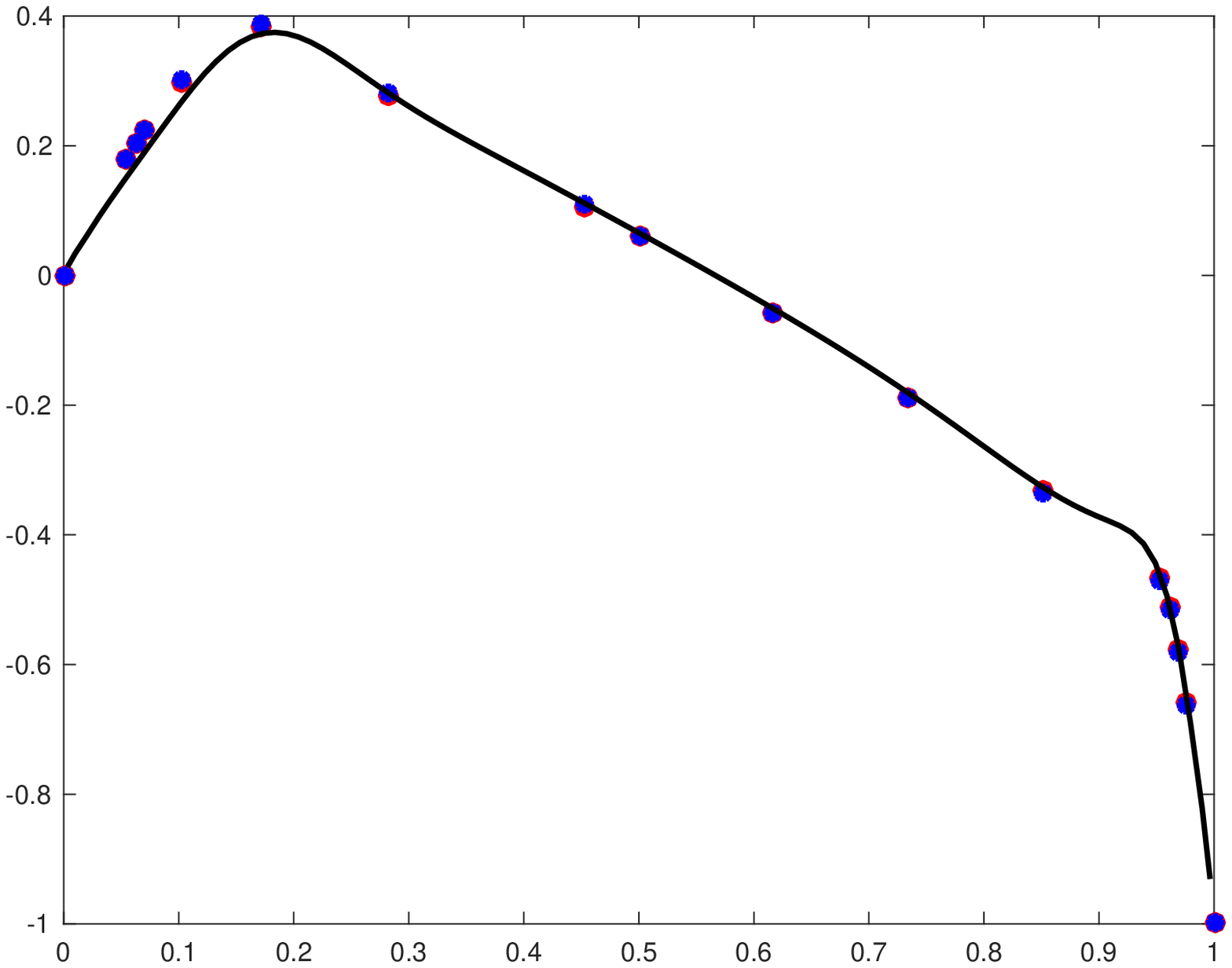}
  \includegraphics[width=0.45\textwidth]{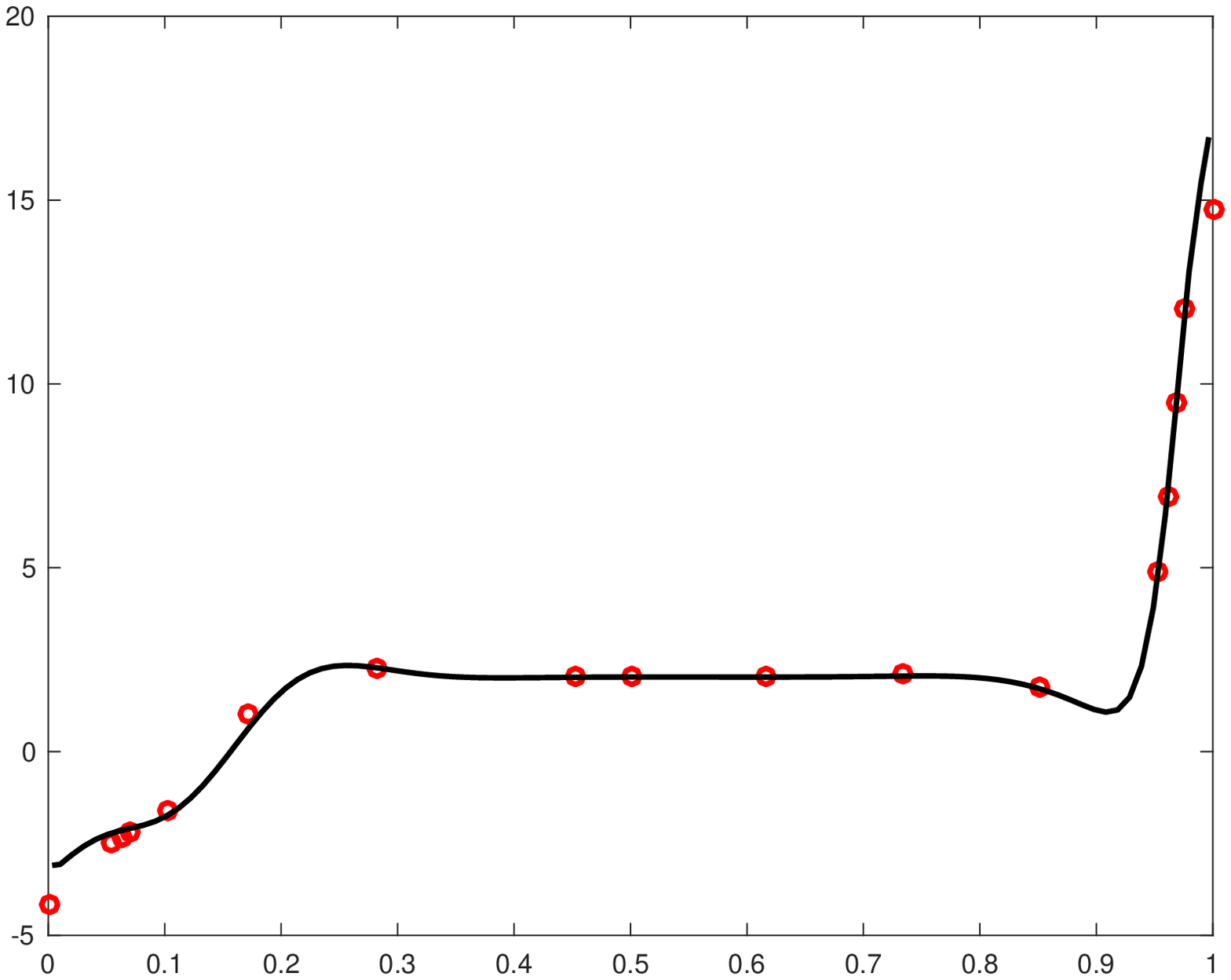}
\caption{Lid driven cavity benchmark at   $Re=1000:$ a)  comparison
with $u$ velocity component values of reference solution b) comparison with $\omega $ values  of reference solutions}
\label{fig:profiles_re1000}
\end{figure}

Results at $Re=1000$ are shown instead in \ref{fig:flow_re1000}. Again, the streamfunction is plotted with the same isoline values  as suggested in \cite{bruneau:2006} and a good general agreement is found between our results and those reported in \cite{botella:1998}, \cite{bruneau:2006}. At a more quantitative level, we show in Figure \ref{fig:profiles_re1000}
the $u$ component of velocity and  the vorticity values along the axis of the
cavity, as computed by the proposed method (black line) on the same mesh described above
with a time step that yields   $ u\Delta t/\Delta x \approx 6, $  $ \nu\Delta t/(2\Delta x^2)\approx 4.$
The results are compared to those of reference solutions presented in \cite{botella:1998}, shown in
blue dots and red circles for velocity and vorticity, respectively. It can be seen that also good
quantitative agreement is achieved, in spite of the relatively coarse mesh, even though
the accuracy loss at the boundaries discussed in Section \ref{fullysl} is apparent in the results
for the vorticity variable.  The maximum horizontal velocity value along the axis of the cavity
was $u_{max}=0.37487,$ while the maximum and minimum vertical velocity values along the 
horizontal centerline of the cavity
were $v_{max}=0.36034,$ $v_{min}=-0.49989,$  respectively. 
The vorticity value at the center of the cavity was computed as  $\omega_{cen}=2.02641.$ 
These values imply  relative
errors of order $10^{-2}$ with respect to the reference steady state solution in \cite{botella:1998}.

         \section{Conclusions}
 \label{conclu} \indent
 We have presented a fully semi-Lagrangian, explicit time discretization method for the Navier Stokes equations
 in the vorticity streamfunction formulation, which allows to use time steps much longer than
 conventional explicit discretizations without the need to solve any linear systems apart from the one resulting
 from the discretization of the Poisson equation. While only first order in time, the method allows in practice
 to achieve a sufficient accuracy for many applications, since high order
 interpolants can be used to reduce   numerical dissipation in a very straightforward way.
 
The effectiveness of the method was demonstrated by numerical results on classical benchmarks
of two dimensional incompressible flows. In ongoing work, this approach is being
 extended to the pressure--velocity formulation
 of the Navier Stokes equations. Furthermore, we are aiming to implement this approach in the framework of
 the high order adaptive SL-DG discretizations proposed in  \cite{tumolo:2015}, \cite{tumolo:2013}.
   
\section*{Acknowledgements}

Part of this research work has been financially supported by the INDAM--GNCS project
 \textit{``Metodi numerici semi-impliciti e semi-Lagrangiani per sistemi iperbolici di leggi di bilancio''} and by University Roma Tre.

\bibliographystyle{plain}

\bibliography{ns_slag}

\end{document}